\documentclass [11pt]{amsart}
\usepackage{amsmath,amssymb}
 \usepackage{amsthm, amsfonts}
 \usepackage{enumerate}
 \usepackage{amssymb,amsmath,amsthm, amsfonts}
 \usepackage{cite}
\usepackage[pagewise]{lineno}
\newtheorem{theorem}{Theorem}[section]
\newtheorem{proposition}[theorem]{Proposition}
\newtheorem{lemma}[theorem]{Lemma}
\newtheorem{corollary}[theorem]{Corollary}
\theoremstyle{definition}

\newtheorem{example}[theorem]{Example}

\theoremstyle{remark}

\numberwithin{equation}{section}

\begin{document}

\title [{{The Zariski topology on the graded primary spectrum }}]{The Zariski topology on the graded primary spectrum of a graded module over a graded commutative ring }

 \author[{{S. Salam and K. Al-Zoubi,  }}]{\textit{ Saif Salam and Khaldoun Al-Zoubi}*}

\address
{\textit{ Saif Salam, Department of Mathematics and
Statistics, Jordan University of Science and Technology, P.O.Box
3030, Irbid 22110, Jordan.}}
\bigskip
{\email{\textit{smsalam19@sci.just.edu.jo}}}

\address
{\textit{Khaldoun Al-Zoubi, Department of Mathematics and
Statistics, Jordan University of Science and Technology, P.O.Box
3030, Irbid 22110, Jordan.}}
\bigskip
{\email{\textit{kfzoubi@just.edu.jo}}}


\date{}
\begin{abstract}
Let $R$ be a $G$-graded ring and M be a $G$-graded $R$-module. We define the graded primary spectrum of $M$, denoted by $\mathcal{PS}_G(M)$, to be the set of all graded primary submodules $Q$ of M such that $(Gr_M(Q):_R M)=Gr((Q:_R M))$. In this paper, we define a topology on $\mathcal{PS}_G(M)$  having the Zariski topology on the graded prime spectrum $Spec_G(M)$ as a subspace topology, and investigate several topological properties of this topological space.
\end{abstract}

\keywords{graded primary submodules, graded primary spectrum, Zariski topology. \\
$*$ Corresponding author}
 \maketitle


 \section{Introduction and Preliminaries }
Let $G$ be a multiplicative group with identity $e$ and R be a commutative ring with identity. Then R is called a $G$-graded ring if there exist additive subgroups $R_g$ of $R$ indexed by the elements $g\in G$ such that $R=\underset{g\in G}{\oplus R_g} $ and $R_g R_h \subseteq R_{gh}$ for all $g,h\in G$. The elements of $R_g$ are called homogeneous of degree g. If $r\in R$, then $r$ can be written uniquely as $\underset{g\in G}{\sum r_g }$, where $r_g$ is the component of $r$ in $R_g$. The set of all homogeneous elements of R is denoted by $h(R)$, i.e. $h(R)=\underset{g\in G}{\bigcup R_g}$. Let $R$ be a $G$-graded ring and $I$ be an ideal of R. Then I is called $G$-graded ideal of $R$ if $I=\underset{g\in G}{\oplus}{ (I\bigcap R_g)}$. By $I\lhd _G R$,  we mean that $I$ is a $G$-graded ideal of $R$, (see \cite{13}). The graded radical of $I$ is the set of all $a=\underset{g\in G}{\sum}a_g\in R$ such that for each $g\in G$ there exists $n_g > 0$ with $a_g ^{n_g}\in I$. By $Gr(I)$ (resp. $\sqrt{I}$) we mean the graded radical (resp. the radical) of $I$, (see \cite{18}).  The graded prime spectrum $Spec_G(R)$ of a graded ring $R$ consists of all graded prime ideals of $R$. It is known that $Spec_G(R)$ is a topological space whose closed sets are $V_G^R(I)=\lbrace p\in Spec_G(R)$ $\vert$  $I\subseteq p\rbrace$  for each graded ideal $I$ of $R$ (see, for example, \cite{14, 16, 18}).

Let $R$ be a $G$-graded ring and $M$ a left $R$-module. Then $M$ is said to be a $G$-graded $R$-module if $M=\underset{g\in G}{\oplus M_g}$ with $R_g M_h\subseteq M_{gh}$ for all $g,h \in G$, where $M_g$ is an additive subgroup of $M$ for all $g\in G$. The elements of $M_g$ are called homogeneous of degree $g$. If $x\in M$, then $x$ can be written uniquely as $\underset{g\in G}{\sum x_g}$, where $x_g$ is the component of $x$ in $M_g$. The set of all homogeneous elements of $M$ is denoted by $h(M)$, i.e. $h(M)=\underset{g\in G}{\bigcup M_g}$. Let $M=\underset{g\in G}{\oplus M_g}$ be a $G$-graded $R$-module. A submodule $N$ of $M$ is called a $G$-graded $R$-submodule of $M$ if $N=\underset{g\in G}{\oplus}(N\bigcap M_g)$. By $N\leq _G M$ (resp. $N<_G M$) we mean that $N$ is a graded submodule (resp. a proper graded submodule) of $M$, (see \cite{13}). If $N\leq_G M$, then  $(N :_R M)=\{ r\in R\, \mid \, rM\subseteq N \}$  is a graded
ideal of $R$, (see \cite[Lemma 2.1]{3}). A proper graded submodule $P$ of $M$ is called a graded prime submodule of $M$ if whenever $r\in h(R)$ and $m\in h(M)$ with $rm\in P$, then either $m\in P$ or $r\in (P:_R M)$. It is easily seen that, if $P$ is a graded prime submodule of $M$, then $(P:_R M)$ is a graded prime ideal of $R$ (see \cite[Proposition2.7]{3}). The graded prime spectrum of $M$, denoted by $Spec_G(M)$, is the set of all graded prime submodules of $M$.  A proper graded submodule $Q$ of $M$ is called a graded primary submodule of $M$, if  whenever $r\in h(R)$ and $m\in h(M)$ with $rm\in Q$, then either $m\in Q$ or $r\in Gr((Q:_R M))$. Graded prime submodules and Graded primary submodules of graded modules have been studied by various authors (see, for example \cite{1,2,3,4,5,15}). The graded radical of a proper graded submodule $N$ of  $M$, denoted by $Gr_M(N)$, is defined to be the intersection of all graded prime submodules of $M$ containing $N$. If $N$ is not contained in any graded prime submodule of $M$, then $Gr_M(N)=M$, (see\cite{5}).

A graded $R$-module $M$ is called a multiplication graded $R$-module if any $N\leq_G M$ has the form $IM$ for some $I\lhd_G R$. If $N$ is a graded submodule of a multiplication graded module $M$, then $N=(N:_{R}M)M$, (see \cite{15}). A graded submodule $N$ of a graded module $M$ is called graded maximal submodule of $M$ if $N\neq M$ and there is no graded submodule $L$ of $M$ such that $N\subset L\subset M$. A graded $R$-module $M$ is called a cancellation graded module if $IM=JM$ for graded ideals $I$ and $J$ of $R$ implies that $I=J$. A graded ring $R$ is called graded integral domain, if whenever $ab=0$ for $a,b\in h(R)$, then $a=0$ or $b=0$. A graded principal ideal domain $R$ is a graded integral domain in which every graded ideal of $R$ is generated by a homogeneous element. One can easily see that, if $R$ is a graded principal ideal domain, then every non-zero graded prime ideal of $R$ is graded maximal.

Let $M$ be a $G$-graded $R$-module and  let $\zeta ^ * (M)=\{V_G^* (N)\, \mid \, N\leq_G M\}$ where $V_G^*(N)=\{P\in Spec_G(M)\, \mid \, N\subseteq P\}$ for any $N\leq _G M$. Then $M$ is called a $G$-top module if the set $\zeta^*(M)$ is closed under finite union. In this case, $\zeta^*(M)$ generates a topology on $Spec_G (M)$ and this topology is called the quasi Zariski topology on $Spec_G(M)$. In contrast with $\zeta^*(M)$, $\zeta(M)=\{V_G(N)\, \mid \, N\leq_G M\}$ where $V_G(N)=\{P\in Spec_G(M)\, \mid \, (N:_R M)\subseteq(P:_R M)\}$ for any $N\leq_G M$ always generates a topology on $Spec_G(M)$. Let $M$ be a $G$-graded $R$-module. Then the map $\varphi:Spec_G(M)\rightarrow Spec_G(R/Ann(M))$ by $\varphi(P)=(P:_R M)/Ann(M)$ is called the natural map on $Spec_G(M)$. For more details concerning the topologies on $Spec_G(M)$ and the natural map on  $Spec_G(M)$, one can look in \cite{7,14}.

In this paper, we call the set of all graded primary submodules $Q$ of a graded module $M$ satisfying the condition $(Gr_M(Q):_R M)=Gr((Q:_R M))$ the graded primary spectrum of $M$ and denote it by $\mathcal{PS}_G(M)$. It is easy to see that $Spec_G(M)\subseteq \mathcal{PS}_G(M)$. The converse inclusion is not always true. For example, if $F$ is a $G$-graded field and $M$ is a $G$-graded $F$-module, then $Spec_G(M)=\mathcal{PS}_G(M)=\{N\, \mid \, N<_G M\}$. But if we take the ring of integers $R=\mathbb{Z}$ as a $\mathbb{Z}_2$-graded $\mathbb{Z}$-module by $R_0=\mathbb{Z}$ and $R_1=\{0\}$, then $4\mathbb{Z}\in \mathcal{PS}_{\mathbb{Z}_2}(\mathbb{Z})-Spec_{\mathbb{Z}_ 2}(\mathbb{Z})$. For a $G$-graded $R$-module M, it is clear that $Gr_M(Q)\neq M$ for any $Q \in \mathcal{PS}_G(M)$ as $Gr((Q:_R M))\in Spec_G(R)$. We introduce the primary $G$-top module which is a generalization of the $G$-top module. For this, we define the variety of any $N\leq_G M$ by $\nu^*_G(N)=\{Q\in \mathcal{PS}_G(M)\, \mid \, N\subseteq Gr_M(Q)\}$ and we set $\Omega^*(M)=\{\nu^*_G (M)\, \mid \, N\leq_G M\}$. Then $M$ is called a primary $G$-top module if $\Omega^*(M)$ is closed under finite union. When this the case, the topology generated by $\Omega^*(M)$ is called the quasi-Zariski topology on $\mathcal{PS}_G (M)$. In particular, every primary $G$-top module is a $G$-top module. Next, we define another variety of any $N\leq_G M$ by $\nu_G (N)=\{Q\in \mathcal{PS}_G(M)\, \mid \, (N:_R M)\subseteq (Gr_M(Q):_R M)\}$. Then the collection $\Omega(M)=\{\nu_G (N) \, \mid \, N\leq_G N\}$ satisfies the axioms for closed sets of a topology on $\mathcal{PS}_G (M)$, which is called the Zariski topology on $\mathcal{PS}_G (M)$, or simply $\mathcal{PZ}_G$-topology. We give some properties of these topologies. We also relate some properties of the graded primary spectrum $\mathcal{PS}_G (M)$ and $Spec_G(R/Ann(M))$  by introducing the map $\rho :\mathcal{PS}_G (M)\rightarrow Spec_G (R/Ann(M))$ given by $\rho(Q)=(Gr_M(Q):_R M)/Ann(M)$. It should be noted that $(Gr_M(Q):_R M)\in Spec_G (R)$, since $Q$ is a graded primary submodule of $M$ and $(Gr_M(Q):_R M)=Gr((Q:_R M))$. In the last two sections, we find a base for the Zariski topology on $\mathcal{PS}_G(M)$ and we obtain some observations and results that concern some conditions under which $\mathcal{PS}_G (M)$ is quasi-compact, irreducible, $T_0$-space or spectral space.

Throughout this paper, $G$ is a multiplicative group, $R$ is a commutative $G$-graded ring with identity and $M$ is a $G$-graded $R$-module. We assume that $Spec_G (M)$ and $\mathcal{PS}_G (M)$ are non-empty.

\section{The Zariski topology on $\mathcal{PS}_G(M)$}

 In this section, we introduce different varieties for graded submodules of graded modules. Using the properties of these varieties, we define the quasi Zariski topology and the $\mathcal{PZ}_G$-topology on $\mathcal{PS}_G (M)$. We also give some relationships between $\mathcal{PS}_G(M)$, $Spec_G(R/Ann(M))$ and $Spec_G(M)$.

\begin{theorem} Let $M$ be a $G$-graded $R$-module. For any $G$-graded submodule $N$ of $M$, we define the variety of $N$ by $\nu^*_G(N)=\{Q\in \mathcal{PS}_G(M)\, \mid \, N\subseteq Gr_M(Q)\}$. Then the following hold:
 \begin{enumerate}

\item$\nu^*_G(0)=\mathcal{PS}_G(M)$ and $\nu^*_G(M)=\emptyset$.
\item If $N,N^\prime\leq_G M$ and $N\subseteq N^\prime$, then $\nu^*_G(N^\prime)\subseteq \nu^*_G(N)$.
\item $\underset{i\in I}{\bigcap} \nu^*_G(N_i)=\nu^*_G(\underset{i\in I}{\sum}N_i)$ for any indexing set I and any family of graded submodules $\{N_i\}_{i\in I}$.
\item $\nu^*_G (N) \cup \nu^*_G(N^\prime)\subseteq \nu^*(N\cap N^\prime)$ for any $N, N^\prime \leq_G M$.
\item $\nu^*_G(N)=\nu^*_G(Gr_M(N))$ for any $N\leq_G M$.
 \end{enumerate}
\end{theorem}
\begin{proof}
(1) and (2) are obvious. \\
(3) Since $N_i\subseteq \underset{i\in I}{\sum}N_i$ for all $i\in I$, then by (2) we have $\nu^*_G(\underset{i\in I}{\sum}N_i)\subseteq \nu^*_G(N_i)$ for all $i\in I$. Therefore $\nu^*_G(\underset{i\in I}{\sum}N_i)\subseteq \underset{i\in I}{\bigcap} \nu^*_G(N_i)$. Conversely, let $Q \in \underset{i\in I}{\bigcap}\nu^*_G(N_i)$. Then $N_i\subseteq Gr_M(Q)$ for all $i\in I$, which implies that $\underset{i\in I}{\sum}N_i \subseteq Gr_M(Q)$. Hence $Q\in \nu^*_G(\underset{i\in I}{\sum}N_i)$.  \\
(4) Since $N\cap N^\prime \subseteq N$ and $N\cap N^\prime \subseteq N^\prime$, then by (2) we have $\nu^*_G(N)\subseteq \nu^*_G(N\cap N^\prime)$ and $\nu^*_G(N^\prime)\subseteq \nu^*_G(N\cap N^\prime)$. Therefore $\nu^*_G(N)\cup \nu^*_G(N^\prime)\subseteq \nu^*_G(N\cap N^\prime)$. \\
(5) As $N\subseteq Gr_M(N)$, we obtain $\nu^*_G(Gr_M(N))\subseteq \nu^*_G(N)$. Conversely, let $Q\in \nu^*_G(N)$. Then $N\subseteq Gr_M(Q)$. So $Gr_M(N)\subseteq Gr_M(Gr_M(Q))=Gr_M(Q)$. Thus $Q\in \nu^*_G(Gr_M(N))$.
\end{proof}

Note that the reverse inclusion in Theorem 2.1 (4) is not always true. Take $R=\mathbb{Z}, G=\mathbb{Z}_2, M=\mathbb{Z}\times \mathbb{Z}, N=4\mathbb{Z}\times \{0\}$ and $N^\prime=\{0\}\times 4\mathbb{Z}$. Then $R$ is a $G$-graded ring by $R_0=\mathbb{Z}$ and $R_1=\{0\}$. Also $M$ is a $G$-graded $R$-module by $M_0=\mathbb{Z}\times \{0\}$ and $M_1=\{0\}\times \mathbb{Z}$. Moreover, $N, N^\prime \leq_G M$. Now, $\nu^*_G(N\cap N^\prime)=\nu^*_G(\{(0,0)\})=\mathcal{PS}_G(M)$. Let $P=\{(0,0)\}$. Then $P\in Spec_G(M)\subseteq \mathcal{PS}_G(M)$. It follows that $P\in \nu^*_G(N\cap N^\prime)$ and $Gr_M(P)=P$. But $N\nsubseteq P$ and $N^\prime \nsubseteq P$. Thus $P\notin \nu^*_G(N) \cup \nu^*_G(N^\prime)$.

By Theorem 2.1 (1), (3), and (4), the collection $\Omega^*(M)= \{\nu^*_G(N)\, \mid \, N\leq_G M\}$ satisfies the axioms for closed sets of a topology on $\mathcal{PS}_G(M)$ if and only if $\Omega^*(M)$ is closed under finite union. When this is the case, we call $M$ a primary $G$-top module and we call the generated topology the quasi Zariski topology on $\mathcal{PS}_G(M)$, or $\mathcal{PZ}_G^q$-topology for short. It is clear that, every $G$-graded simple $R$-module is a primary $G$-top module. In the following theorem, we show that every multiplication graded $R$-module is a primary $G$-top module.

\begin{theorem}
If $M$ is a multiplication graded $R$-module, then $M$ is a primary $G$-top module.
\end{theorem}
\begin{proof}
Let $N, N^\prime \leq_G M$. It is sufficient to show that $\nu^*_G(N\cap N^\prime)\subseteq \nu^*_G(N)\cup \nu^*_G(N^\prime)$. So let $Q\in \nu^*_G(N\cap N^\prime )$. Then $N\cap N^\prime \subseteq Gr_M(Q)$. It follows that $(N:_{R}M)\cap (N^\prime:_R M)=(N\cap N^\prime:_R M)\subseteq (Gr_M(Q):_R M)=Gr((Q:_R M))\in Spec_G(R)$ as $Q$  is a graded primary submodule. Therefore $(N:_R M)\subseteq (Gr_M(Q):_R M)$ or $(N^\prime:_R M)\subseteq (Gr_M(Q):_R M)$. Since $M$ is a multiplication graded $R$-module, then $N=(N:_R M)M\subseteq (Gr_M(Q):_R M)M=Gr_M(Q)$ or $N^\prime=(N^\prime:_R M)M\subseteq (Gr_M(Q):_R M)M=Gr_M(Q)$. Thus $Q\in \nu^*_G(N)\cup \nu^*_G(N^\prime)$.
\end{proof}

\begin{proposition}
 Let $M$ be multiplication $G$-graded $R$-module, $N\leq_G M$ and  $I,J\lhd_G R$. Then the following hold:
  \begin{enumerate}
\item $\nu^*_G(N)\cup \nu^*_G(IM)= \nu^*_G(IN)$.
\item $\nu^*_G(IM) \cup \nu^*_G(JM)=\nu^*_G((IJ)M)$.
\end{enumerate}
 \end{proposition}

\begin{proof} (1) Since $IN\subseteq N$ and $IN\subseteq IM$, then $\nu^*_G(N)\subseteq \nu^*_G(IN)$ and $\nu^*_G(IM)\subseteq \nu^*_G(IN)$. Thus $\nu^*_G(N)\cup \nu^*_G(IM)\subseteq \nu^*_G(IN)$. To establish the reverse inclusion, let $Q\in \nu^*_G(IN)$ such that $IN\subseteq Gr_M(Q)$. As $Q\in {PS}_G(M)$ and $M$ is a multiplication graded module, then by \cite[Theorem 14]{15}, $Gr_M(Q)\in Spec_G(M)$. It follows that $N\subseteq Gr_M(Q)$ or $I\subseteq (Gr_M(Q):_R M)$. Therefore $N\subseteq Gr_M(Q)$ or $IM\subseteq Gr_M(Q)$ which completes the proof.\\
(2) It is clear by (1).
\end{proof}

Now we define another variety for a graded submodule $N$ of a graded $R$-module $M$. We set $\nu_G(N)=\{Q\in {PS}_G(M)\, \mid\, (N:_R M)\subseteq (Gr_M(Q):_R M)\}$. We state some properties of this variety in the following theorem to construct the Zariski topology on $\mathcal{PS}_G(M)$.

\begin{theorem} Let $M$ be a $G$-graded $R$-module and let $N, N^\prime, N_i\leq_G M$ for any $i\in I$, where $I$ is an indexing set. Then the following hold:
 \begin{enumerate}
\item $\nu_G(0)=\mathcal{PS}_G(M)$ and $\nu_G(M)=\emptyset$.
\item $\underset{i\in I}{\bigcap}\nu_G(N_i)=\nu_G(\underset{i\in I}{\sum}(N_i:_R M)M)$.
\item $\nu_G(N)\cup \nu_G(N^\prime)=\nu_G(N\cap N^\prime)$.
\item If $N\subseteq N^\prime$, then $\nu_G(N^\prime)\subseteq\nu_G(N)$.
\end{enumerate}
\end{theorem}
\begin{proof}
(1) and (4) are trivial.\\
(2) Let $Q\in \underset{i\in I}{\bigcap}\nu_G(N_i)$. Then $(N_i:_R M)\subseteq (Gr_M(Q):_R M)$ for all $i\in I$. So $(N_i :_R M)M\subseteq (Gr_M(Q):_R M)M$ for all $i\in I$. This implies that $\underset{i\in I}{\sum}(N_i:_R M)M \subseteq(Gr_M(Q):_R M)M$. Therefore $(\underset{i\in I}{\sum}(N_i:_R M)M : M)\subseteq((Gr_M(Q):_R M)M :_R M)=(Gr_M(Q):_R M)$. Hence $Q\in \nu_G(\underset{i\in I}{\sum}(N_i:_R M)M)$. For the reverse inclusion, let $Q\in \nu_G(\underset{i\in I}{\sum}(N_i:_R M)M)$. Then $(\underset{i\in I}{\sum}(N_i :_R M)M:_R M)\subseteq (Gr_M(Q):_R M)$. But for any $i\in I$, we have $(N_i :_R M)=((N_i:_R M)M:_R M)\subseteq (\underset{i\in I}{\sum}(N_i:_R M)M:_R M)\subseteq (Gr_M(Q):_R M)$. Thus $Q\in \underset{i\in I}{\bigcap}\nu_G(N_i)$. \\
(3) For any $Q\in \mathcal{PS}_G(M)$, we have $Q\in \nu_G(N\cap N^\prime)\Leftrightarrow (N\cap N^\prime:_R M)\subseteq (Gr_M(Q):_R M)=Gr((Q:_R M))\Leftrightarrow (N:_R M)\cap(N^\prime:_R M)\subseteq Gr((Q:_R M))\in Spec_G(R)\Leftrightarrow (N:_R M)\subseteq Gr((Q:_R M))$ or $(N^\prime:_R M)\subseteq Gr((Q:_R M))\Leftrightarrow Q\in \nu_G(N)\cup \nu_G(N^\prime)$. Hence $\nu_G(N\cap N^\prime)=\nu_G(N)\cup \nu_G(N^\prime)$.
\end{proof}

In view of Theorem 2.4 (1), (2) and (3), the collection $\Omega(M)=\{\nu_G(N)\, \mid \, N\leq_G M\}$ satisfies the axioms for closed sets of a topology on $\mathcal{PS}_G(M)$, which is called the primary Zariski topology on $\mathcal{PS}_G(M)$, or $\mathcal{PZ}_G$-topology for short. In the next proposition, we set some conditions under which $\nu_G(N)=\nu_G(Gr_M(N))$ for $N\leq_G M$.

\begin{proposition} Let $M$ be a $G$-graded $R$-module and $N\leq_G M$. If $N\in \mathcal{PS}_G(M)$ or $M$ is multiplication graded module, then $\nu_G(N)=\nu_G (Gr_M(N))$.
\end{proposition}
\begin{proof}
It is clear that, $\nu_G(Gr_M(N))\subseteq \nu_G(N)$. So let $Q\in \nu_G(N)$ and it is enough to show that $(Gr_M(N):_R M)\subseteq (Gr_M(Q):_R M)$. Firstly, assume that $ N\in \mathcal{PS}_G(M)$. Since $Q\in \nu_G(N)$, then $(N:_R M)\subseteq Gr((Q:_R M))$. As $N\in \mathcal{PS}_G(M)$, we obtain $(Gr_M(N):_R M)=Gr((N:_R M))\subseteq Gr(Gr((Q:_R M)))=Gr((Q:_R M))=(Gr_M(Q):_R M)$ which completes the proof of the first case. For the second case, assume that $M$ is multiplication graded module. Again, since $Q\in \nu_G(N)$, then $(N:_R M)\subseteq (Gr_M(Q):_R M)$. This follows that $N=(N:_R M)M\subseteq (Gr_M(Q):_R M)M=Gr_M(Q)$, which implies that $Gr_M(N)\subseteq Gr_M(Gr_M(Q))=Gr_M(Q)$. Hence $(Gr_M(N):_R M)\subseteq (Gr_M(Q):_R M)$.
\end{proof}
Now we state some relations between the varieties $V^*_G(N)$, $V_G(N)$, $\nu^*_G(N)$ and $\nu_G(N)$ for any graded submodule $N$ of an $R$-module $M$. These relations will be used continuously throughout the rest of this paper.\\
\begin{lemma}
Suppose that $N$ and $N^\prime$ are graded submodules of a $G$-graded $R$-module $M$ and that $I$ is a $G$-graded ideal of $R$. Then the following hold:
\end{lemma}
\begin{enumerate}
\item $V_G(N)=\nu_G(N)\cap Spec_G(M)$.
\item $V^*_G(N)=\nu^*_G(N)\cap Spec_G(M)$.
\item If $Gr((N:_R M))=Gr((N^\prime:_R M))$, then $\nu_G(N)=\nu_G(N^\prime)$. The converse is also true if $N, N^\prime\in \mathcal{PS}_G(M)$.
\item $\nu_G(N)=\nu_G((N:_R M)M)=\nu^*_G((N:_R M)M)=\nu^*_G(Gr((N:_R M))M)$. In particular, $\nu^*_G(IM)=\nu_G(IM)$.
\end{enumerate}
\begin{proof}
The proof is straightforward.
\end{proof}
\begin{corollary}
Every primary $G$-top module is a $G$-top module.
\end{corollary}
\begin{proof}
 Let $M$ be a primary $G$-top module and $N, N^\prime \leq_G M$. By Lemma 2.6 (2), we have $V^*_G(N)\cup V^*_G(N^\prime)= (Spec_G(M) \cap \nu^*_G(N))\cup (Spec_G(M)\cap \nu^*_G(N^\prime))=Spec_G(M)\cap(\nu^*_G(N)\cup \nu^*_G(N^\prime))=Spec_G(M)\cap \nu^*_G(J)=V^*_G(J)$ for some graded submodule $J$ of $M$ and hence $M$ is a $G$-top module.
 \end{proof}

By Corollary 2.7, if $M$ is a primary $G$-top module, then $\zeta^*(M)=\{V^*_G(N)\, \mid \, N\leq_G M\}$ induces the quasi Zariski topology on $Spec_G(M)$ which will be, by Lemma 2.6 (2), a topological subspace of $\mathcal{PS}_G(M)$ equipped with $\mathcal{PZ}_G^q$-topology. Also by Lemma 2.6 (1), $Spec_G(M)$ with the Zariski topology is a topological subspace of $\mathcal{PS}_G(M)$ equipped with the $\mathcal{PZ}_G$-topology for any G-graded $R$-module $M$.

Consider $\varphi$ and $\rho$ as described in the introduction. Let $M$ be a $G$-graded $R$-module. For $p\in Spec_G(R)$, we set $\mathcal{PS}_G^p(M)=\{Q\in \mathcal{PS}_G(M) \, \mid \, (Gr_M(Q):_R M)=p\}$.

\begin{proposition}
 The following statements are equivalent for any $G$-graded $R$-module $M$:
  \begin{enumerate}
\item If whenever $Q, Q^\prime \in PS_G(M) $ with $\nu_G(Q)=\nu_G(Q^\prime)$, then $Q=Q^\prime$.
\item $|\mathcal{PS}_G^p(M)|\leq 1$ for every $p\in Spec_G(R)$.
\item $\rho$ is injective.
\end{enumerate}
\end{proposition}
\begin{proof}(1)$\Rightarrow$(2): Let $p\in Spec_G(R)$ and  $Q, Q^\prime \in \mathcal{PS}_G^q(M)$. Then $Q, Q^\prime \in \mathcal{PS}_G(M)$ and $(Gr_M(Q):_R M)=(Gr_M(Q^\prime):_R M)=p$. By Lemma 2.6 (3), we have $\nu_G(Q)=\nu_G(Q^\prime)$. So, by the assumption (1), $Q=Q^\prime$.\\
(2)$\Rightarrow$(3): Assume that $\rho(Q)=\rho(Q^\prime)$, where $Q, Q^\prime\in \mathcal{PS}_G(M)$. Then $(Gr_M(Q):_R M)=(Gr_M(Q^\prime):_R M)$. Let $p=(Gr_M(Q):_R M)\in Spec_G(R)$. Then we get $Q, Q^\prime \in \mathcal{PS}_G^p(M)$ and by the hypothesis we obtain $Q=Q^\prime$.\\
(3)$\Rightarrow$(1): Let $Q, Q^\prime \in \mathcal{PS}_G(M)$ with $\nu_G(Q)=\nu_G(Q^\prime)$. Then, by Lemma 2.6 (3), we have $(Gr_M(Q):_R M)=(Gr_M(Q^\prime):_R M)$. So $\rho(Q)=\rho(Q^\prime)$. Since $\rho$ is injective, then $Q=Q^\prime$.
\end{proof}
\begin{corollary}
If $|\mathcal{PS}_G^p(M)|=1$ for every $p\in Spec_G(R)$, then $\rho$ is bijective.
\end{corollary}
\begin{proof} It is clear by Proposition 2.8.
\end{proof}
 Let $M$ be a $G$-graded $R$-module. From now on, we will denote $R/Ann(M)$ by  $\overline{R}$ and any graded ideal $I/Ann(M)$ of $\overline{R}$ by $\overline{I}$. By \cite[Proposition 3.13]{14}, the natural map $\varphi$ of $Spec_G(M)$ is continuous and $\varphi ^{-1}(V_{G}^{R}(\overline{I}))=V_{G}(IM)$, for every graded ideal $I$ of $R$ containing $Ann(M)$. Also, in \cite[Proposition 3.15]{14}, it has been shown that if $\varphi $ is surjective, then $\varphi $ is both open and closed with $\varphi (V_{G}(N))=V_{G}{^{\overline{R}}}(\overline{(N:_{R}M)})$ and $\varphi
(Spec_{G}(M)-V_{G}(N))=Spec_{G}(\overline{R})-V_{G}{^{\overline{R}}}(%
\overline{(N:_{R}M)})$, for any $N\leq _{G}M$. These two results are important in the rest of this section. In the next two propositions, we give
similar results for $\rho$.

\begin{proposition}
 Let $M$ be a $G$-graded $R$-module. Then $\rho^{-1}(V_G^{\overline{R}}(\overline{I}))=\nu_G(IM)$, for every graded ideal $I$ of $R$ containing $Ann(M)$. Therefore $\rho$ is continuous mapping.
 \end{proposition}
\begin{proof}
For any $Q\in \mathcal{PS}_G(M)$, we have $Q\in \rho^{-1}(V_G^{\overline{R}}(\overline{I}))\Leftrightarrow \rho(Q)\in V_G^{\overline{R}}(\overline{I})\Leftrightarrow \overline{I}\subseteq \overline{(Gr_M(Q):_R M)}\Leftrightarrow I\subseteq (Gr_M(Q):_R M)\Leftrightarrow IM\subseteq (Gr_M(Q):_R M)M\Leftrightarrow (IM:_R M)\subseteq ((Gr_M(Q):_R M)M:_R M)=(Gr_M(Q):_R M)\Leftrightarrow Q\in \nu_G(IM)$. Hence $\rho^{-1}(V_G^{\overline{R}}(\overline{I}))=\nu_G(IM)$.
\end{proof}
\begin{proposition}
Let $M$  be a $G$-graded $R$-module. If $\rho$ is surjective, then $\rho$ is both open and closed; more precisely, for any $N\leq_G M$, $\rho(\nu_G(N))=V_G^{\overline{R}}(\overline{(N:_R M)})$ and $\rho(\mathcal{PS}_G(M) - \nu_G(N))=Spec_G(\overline{R}) - V_G^{\overline{R}}(\overline{(N:_R M)})$.
\end{proposition}
 \begin{proof}
By Proposition 2.10, we have $\rho^{-1}(V_G^{\overline{R}}(\overline{I}))=\nu_G(IM)$ for every $I\lhd_G R$ containing $Ann(M)$. So $\rho^{-1}(V_G^{\overline{R}}(\overline{(N:_R M)}))=\nu_G((N:_R M)M)=\nu_G(N)$ for any $N\leq_G M$. It follows that $V_G^{\overline{R}}(\overline{(N:_R M)})=\rho(\rho^{-1}(V_G^{\overline{R}}(\overline{(N:_R M)})))=\rho(\nu_G(N))$ as $\rho$ is surjective. For the second part, note that $\mathcal{PS}_G(M)-\nu_G(N)=\mathcal{PS}_G(M)-\rho^{-1}(V_G^{\overline{R}}(\overline{(N:_R M)}))=\rho^{-1}(Spec_G(\overline{R}))-\rho^{-1}(V_G^{\overline{R}}(\overline{(N:_R M)}))=\rho^{-1}(Spec_G(\overline{R})-V_G^{\overline{R}}(\overline{(N:_R M)}))$. This implies that $\rho(\mathcal{PS}_G(M)-\nu_G(N))=\rho(\rho^{-1}(Spec_G(\overline{R})-V_G^{\overline{R}}(\overline{(N:_R M)})))=Spec_G(\overline{R})-V_G^{\overline{R}}(\overline{(N:_R M)})$.
  \end{proof}
\begin{corollary} Let $M$ be a $G$-graded $R$-module. Then $\rho$ is bijective if and only if $\rho$ is a homeomorphism.
\end{corollary}
\begin{theorem}
Let $M$ be a $G$-graded $R$-module. Consider the following statements:
\begin{enumerate}
\item $Spec_G(M)$ is connected.
\item $\mathcal{PS}_G(M)$ is connected.
\item $Spec_G(\overline{R})$ is connected.
\end{enumerate}
\begin{enumerate}[\upshape (i)]
\item If $\rho$ is surjective, then (1)$\Rightarrow$(2)$\Leftrightarrow$(3).
 \item If $\varphi$ is surjective, then all the three statements are equivalent.
 \end{enumerate}
\end{theorem}
\begin{proof}
(i) (1)$\Rightarrow$(2): Assume that $Spec_G(M)$ is connected. If $\mathcal{PS}_G(M)$ is disconnected, then there is a $U$ clopen in $\mathcal{PS}_G(M)$ such that $U\neq \emptyset$ and $U\neq \mathcal{PS}_G(M)$. Since $U$ is clopen in $\mathcal{PS}_G(M)$, then $U= \mathcal{PS}_G(M)-\nu_G(N_1)=\nu_G(N_2)$ for some $N_1, N_2\leq_G M$. By Proposition 2.11, we have $\rho(U)$ is clopen in $Spec_G(\overline{R})$. But $\varphi$ is continuous. So $\varphi^{-1}(\rho(U))$ is clopen in $\mathcal{PS}_G(M)$, and so $\varphi^{-1}(\rho(U))=\emptyset$ or $\varphi^{-1}(\rho(U))=Spec_G(M)$. If $\varphi^{-1}(\rho(U))=\emptyset$, then $\varphi^{-1}(\rho(\nu_G(N_2)))=\emptyset$. Thus $\varphi^{-1}(V_G^{\overline{R}}(\overline{(N_2:_R M)}))=\emptyset$. It follows that $V_G(N_2)=\emptyset$, which means that $(N_2:_R M)\nsubseteq(P:_R M)$ for any $P\in Spec_G(M)$. As $U=\nu_G(N_2)\neq \emptyset$, then $\exists Q \in \mathcal{PS}_G(M)$ such that $(N_2:_R M)\subseteq(Gr_M(Q):_R M)$. Since $Gr_M(Q)\neq M$, then $\exists P^\prime \in Spec_G(M)$  such that $Q\subseteq P^\prime$. Therefore $(Gr_M(Q):_R M)\subseteq (Gr_M(P^\prime):_R M)=(P^\prime:_R M)$. Hence $(N_2:_R M)\subseteq (P^\prime:_R M)$ which is a contradiction. Now, if $\varphi^{-1}(\rho(U))=Spec_G(M)$, then $Spec_G(M)=\varphi^{-1}(\rho(\nu_G(N_2)))=V_G(N_2)=Spec_G(M)\cap \nu_G(N_2)$. It follows that $Spec_G(M)\subseteq \nu_G(N_2)=U=\mathcal{PS}_G(M)-\nu_G(N_1)$. As $U=\mathcal{PS}_G(M)-\nu_G(N_1)\neq \mathcal{PS}_G(M)$, then $\exists Q\in \mathcal{PS}_G(M) \cap \nu_G(N_1)$. Therefore $(N_1:_R M)\subseteq (Gr_M(Q):_R M)$ and $\exists P\in Spec_G(M)$ such that $Q\subseteq P$. It follows that $(N_1:_R M)\subseteq (Gr_M(Q):_R M)\subseteq (Gr_M(P):_R M)$. Then $P\in \nu_G(N_1)$. But $P\in Spec_G(M)\subseteq \mathcal{PS}_G(M)-\nu_G(N_1)$. Therefore $P\notin \nu_G(N_1)$ which is also a contradiction. Consequently, $\mathcal{PS}_G(M)$ is a connected space.\\
(2)$\Rightarrow$(3): Since $\rho$ is continuous surjective map and $\mathcal{PS}_G(M)$ is connected, then $Spec_G(\overline{R})$ is connected.\\
(3)$\Rightarrow$(2): Assume by way of contradiction that $\mathcal{PS}_G(M)$ is disconnected. Then there is a $U$ clopen in $\mathcal{PS}_G(M)$ such that $U\neq \emptyset$ and $U\neq \mathcal{PS}_G(M)$. Since $\rho$ is surjective, then $\rho(U)$ is clopen in $Spec_G(\overline{R})$, and so $\rho(U)=\emptyset$ or $\rho(U)=Spec_G(\overline{R})$ as $Spec_G(\overline{R})$ is connected. Also, since $U$ is open in $\mathcal{PS}_G(M)$, then $U=\mathcal{PS}_G(M)-\nu_G(N)$ for some $N\leq_G M$. Now, if $\rho(U)=Spec_G(\overline{R})$, then $Spec_G(\overline{R})=\rho(\mathcal{PS}_G(M)-\nu_G(N))=Spec_G(\overline{R})-V_G^{\overline{R}}(\overline{(N:_R M)})$ by Proposition 2.11. Thus $V_G^{\overline{R}}(\overline{(N:_R M)})=\emptyset$ which implies that $\emptyset=\rho^{-1}(\emptyset)=\rho^{-1}(V_G^{\overline{R}}(\overline{(N:_R M)}))=\nu_G(N)$. Therefore $\nu_G(N)=\emptyset$ and hence $U=\mathcal{PS}_G(M)$, a contradiction. Also if $\rho(U)=\emptyset$, then $U\subseteq\rho^{-1}(\rho(U))=\emptyset$. It follows that $U=\emptyset$ which is also a contradiction. Therefore $\mathcal{PS}_G(M)$ is connected.\\
(ii) If $\varphi$ is surjective, then it is clear that $\rho$ is surjective and hence (1)$\Rightarrow$(2)$\Leftrightarrow$(3) by (i). So it is enough to show that (3)$\Rightarrow$(1. We prove it in a similar way to proof of (3)$\Rightarrow$(2) using the properties of $\varphi$. So again, we assume by  way of contradiction that $Spec_G(M)$ is disconnected, which means that there is a $U$ clopen in $Spec_G(M)$ such that $U\neq \emptyset$ and $U\neq Spec_G(M)$. As $Spec_G(\overline{R})$ is connected and $\varphi$ is surjective, then $\varphi(U)=\emptyset$ or $\varphi(U)=Spec_G(\overline{R})$. Also the open set $U$ can be written as $U=Spec_G(M)-V_G(N)$ for some $N\leq_G M$. It is clear that if $\varphi(U)=\emptyset$, then $U=\emptyset$ and we have a contradiction. If $\varphi(U)=Spec_G(\overline{R})$, then $Spec_G(\overline{R})=\varphi(Spec_G(M)-V_G(N))=Spec_G(\overline{R})-V_G^{\overline{R}}(\overline{(N:_R M)})$, and so $V_G^{\overline{R}}(\overline{(N:_R M)})=\emptyset$. Thus $V_G(N)=\varphi^{-1}(V_G^{\overline{R}}(\overline{(N:_R M)}))=\varphi^{-1}(\emptyset)=\emptyset$. It follows that $U=Spec_G(M)$ which is a contradiction. Therefore $Spec_G(M)$ is a connected space and this completes the proof.
\end{proof}
Let $M$ and $S$ be two $G$-graded $R$-modules. Recall that an $R$-module homomorphism $f: M\rightarrow S$ is called a $G$-graded $R$-module homomorphism if $f(M_g)\subseteq S_g$ for all $g\in G$, see\cite{13}. Let $f:M\rightarrow M^\prime $ be a $G$-graded module epimorphism between the two graded modules $M$ and $M^\prime$. If $N^\prime\leq_G M^\prime$, then $(N^\prime:_R M)=(f^{-1}(N^\prime):_R M)$. Also it is easy to check that, $(N:_R M)=(f(N):_R M)$ and $f(Gr_M(N))=Gr_{M^{\prime}}(f(N))$ for any $N\leq_G M$ containing the kernel of $f$. We will denote the kernel of $f$ by $ker f$.

\begin{lemma}
Let $M$ and $M^\prime$ be $G$-graded $R$-modules. Let $f:M\rightarrow M^\prime$ be a $G$-graded module epimorphism. Then the following hold:
\begin{enumerate}
\item If $Q^\prime \in \mathcal{PS}_G(M^\prime)$, then $f^{-1}(Q^\prime)\in \mathcal{PS}_G(M)$.
\item If $Q\in \mathcal{PS}_G(M)$ and $ker f\subseteq Q$, then $f(Q)\in \mathcal{PS}_G(M^\prime)$.
\end{enumerate}
\end{lemma}
\begin{proof} (1) It is easy to verify that $f^{-1}(Q^\prime)$ is graded primary submodule of $M$ and it remains to show that $(Gr_M(f^{-1}(Q^\prime)):_R M)=Gr((f^{-1}(Q^\prime):_R M))$. Since $ker f\subseteq f^{-1}(Q^\prime)\subseteq Gr_M(f^{-1}(Q^\prime))$ and $(Gr_{M^{\prime}}(Q^\prime):_R M^\prime)=Gr((Q^\prime:_R M^\prime))$, we obtain $(Gr_M(f^{-1}(Q^\prime)):_R M)=(f(Gr_M(f^{-1}(Q^\prime))):_R M^\prime)=(Gr_{M^\prime}(f(f^{-1}(Q^\prime))):_R M^\prime)=(Gr_{M^{\prime}}(Q^\prime):_R M^\prime)=Gr((Q^\prime:_R M^\prime))=Gr((f^{-1}(Q^\prime):_R M))$ as required. \\
(2) First note that $f(Q)$ is a graded proper submodule of $M^\prime$, since $Q$ is a graded proper submodule of $M$ containing $kerf$. Let $rm^\prime \in f(Q)$ for $r\in h(R)$ and $m^\prime \in h(M^\prime)$. As $f$ is a graded module epimorphism and $m^\prime \in h(M^\prime)$, we get $\exists m\in h(M)$ such that $f(m)=m^\prime$, which implies that $f(rm)\in f(Q)$. Thus $\exists t\in Q$ such that $rm-t\in ker f\subseteq Q$. So $rm\in Q$, and so $m\in Q$ or $r\in \sqrt{(Q:_R M)}=\sqrt{(f(Q):_R M^\prime)}$ . Hence $f(Q)$ is a graded primary submodule of $M^\prime$. Moreover, $(Gr_{M^\prime}(f(Q)):_R M^\prime)=(f(Gr_M(Q)):_R M^\prime)=(Gr_M(Q):_R M)=Gr((Q:_R M))=Gr((f(Q):_R M^\prime))$. Therefore $f(Q)\in \mathcal{PS}_G (M^\prime)$.
\end{proof}

\begin{theorem}
Let $M$ and $M^\prime$ be $G$-graded $R$-modules and $f:M\rightarrow M^\prime$ be a graded module epimorphism. Then the mapping $\pi:\mathcal{PS}_G(M^\prime)\rightarrow \mathcal{PS}_G(M)$ by $\pi(Q^\prime)=f^{-1}(Q^\prime)$ is an injective continuous map. Moreover, if $\pi$ is surjective map, then $\mathcal{PS}_G(M)$ is homeomorphic to $\mathcal{PS}_G(M^\prime)$.
\end{theorem}
\begin{proof} By Lemma 2.14, $\pi$ is well-defined. Also, the injectivity of $\pi$ is obvious. Now for any $O\in \mathcal{PS}_G(M^\prime)$ and any closed set $\nu_G(N)$ in $\mathcal{PS}_G(M)$, where $N\leq_G M$, we have $O\in \pi^{-1}(\nu_G(N))=\pi^{-1}(\nu^*_G(Gr((N:_R M))M)) \Leftrightarrow Gr((N:_R M))M\subseteq Gr_M(f^{-1}(O))\Leftrightarrow Gr((N:_R M))\subseteq (Gr_M(f^{-1}(O)):_R M)=Gr((f^{-1}(O):_R M))=Gr((O:_R M^{\prime}))=(Gr_{M^\prime}(O):_R M^\prime)\Leftrightarrow Gr((N:_R M))M^\prime\subseteq Gr_{M^\prime}(O) \Leftrightarrow O\in \nu^*_G(Gr((N:_R M))M^\prime)=\nu_G(Gr((N:_R M))M^\prime)$. Therefore $\pi^{-1}(\nu_G(N))=\nu_G(Gr((N:_R M))M^\prime)$ and hence $\pi$ is continuous. For the last statement, we assume that $\pi$ is surjective and it is enough to show that $\pi$ is closed. So let $\nu_G(N^\prime)$ be a closed set in $\mathcal{PS}_G(M^\prime)$, where $N^\prime \leq_G M^\prime$. As we have seen, $\pi^{-1}(\nu_G(N))=\nu_G(Gr((N:_R M))M^\prime)$ for any $N\leq_G M$. It follows that $\pi^{-1}(\nu_G(f^{-1}(N^\prime)))=\nu_G(Gr((f^{-1}(N^\prime):_R M))M^\prime)=\nu_G(Gr((N^\prime:_R M^\prime))M^\prime)=\nu_G(N^\prime)$ and hence $\pi^{-1}(\nu_G(f^{-1}(N^\prime)))=\nu_G(N^\prime)$. Thus $\pi(\nu_G(N^\prime))= \nu_G(f^{-1}(N^\prime))$ as $\pi$ is surjective. Therefore $\pi$ is closed, and so $\mathcal{PS}_G(M)$ is homeomorphic to $\mathcal{PS}_G(M^\prime)$.
\end{proof}
\begin{corollary}
 Let $M$ and $M^\prime$ be $G$-graded $R$-modules. Let $f:M\rightarrow M^\prime$ be a $G$-graded module isomorphism. Then $\mathcal{PS}_G(M)$ is homeomorphic to $\mathcal{PS}_G(M^\prime)$.
\end{corollary}
\begin{proof}
By Lemma 2.14 (2) and Theorem 2.15.
\end{proof}

\begin{theorem}
Let $R$ be a graded principal ideal domain. Let $M$ be a cancellation multiplication graded $R$-module and $N\leq_G M$. Then $N\in \mathcal{PS}_G(M)$ if and only if $Gr_M(N)\in Spec_G(M)$.
 \end{theorem}

\begin{proof}
Note that if $N\in \mathcal{PS}_G(M)$, then  $Gr_M(N)\in Spec_G(M)$ by \cite[Theorem 14]{15}. Now, suppose that $Gr_M(N)\in Spec_G(M)$. Since $M$ is multiplication graded $R$-module, then $Gr_M(N)=(Gr_M(N):_R M)M$. Also by \cite[Theorem 9]{15}, we have $Gr_M(N)=Gr((N:_R M))M$, which follows that $Gr((N:_R M))M=(Gr_M(N):_R M)M$. Since $M$ is cancellation graded module, then $(Gr_M(N):_R M)=Gr((N:_R M))$ and it remains to prove that $N$ is a graded primary submodule of $M$. If $Gr((N:_R M))=0_R$, then $Gr_M(N)=Gr((N:_R M))M=0_M\in Spec_G(M)$. Thus $N=0_M\in Spec_G(M)$ and hence $N$ is a graded primary submodule of $M$. If $Gr((N:_R M))\neq 0_R$, then $Gr((N:_R M))$ is a graded maximal ideal of $R$ as $R$ is a graded principal ideal domain and $Gr((N:_R M))=(Gr_M(N):_R M)\in Spec_G(R)$. Therefore $(N:_R M)$ is a graded primary ideal of $R$ by \cite[Proposition 1.11]{17}. It follows that $N$ is a graded primary submodule of $M$, by \cite[Theorem 11]{15}, since $M$ is a multiplication graded module.
\end{proof}

\section{A base for the Zariski topology on $\mathcal{PS}_G(M)$}
Let $M$ be a $G$-graded $R$-module. In \cite[Theorem 2.3]{14}, it has been proved that for each $r \in h(R)$, the set $D_r=Spec_G(R)-V_G^R(rR)$ is open in $Spec_G(R)$ and the family $\{D_r\, \mid \, r\in h(R)\}$ is a base for the Zariski topology on $Spec_G(R)$. In addition, each $D_r$ is quasi-compact and thus $D_1=Spec_G(R)$ is quasi-compact. In this section, we set $S_r=\mathcal{PS}_G(M)-\nu_G(rM)$ for each $r\in h(R)$ and prove that $S=\{S_r\,\mid\,r\in h(R)\}$ forms a base for the Zariski-topology on $\mathcal{PS}_G(M)$. Also, we show that each $S_r$ is quasi-compact and hence $\mathcal{PS}_G(M)$ is quasi-compact.

\begin{proposition}
For any $G$-graded $R$-module $M$, the set $S=\{S_r\,\mid\,r\in h(R)\}$ forms a base for the Zariski topology on $\mathcal{PS}_G(M)$.
\end{proposition}
\begin{proof}
 Let $U=\mathcal{PS}_G(M)-\nu_G(N)$ be an open set in $\mathcal{PS}_G(M)$, where $N\leq_G M$. Let $Q\in U$ and it is enough to find an element $r\in h(R)$ such that $Q\in S_r\subseteq U$. Since $Q\in U$, then $(N:_R M)\nsubseteq(Gr_M(Q):_R M)$, and so there exists $x\in R$ and $g\in G$  such that $x_g\in (N:_R M)-(Gr_M(Q):_R M)$. Take $r=x_g\in h(R)$. Therefore $(rM:_R M)\nsubseteq (Gr_M(Q):_R M)$ and thus $Q\in S_r$. Now for any $Q^\prime \in S_r$, we have $(rM:_R M)\nsubseteq (Gr_M(Q^\prime):_R M)$, which implies that $(N:_R M)\nsubseteq (Gr_M(Q^\prime):_R M)$. Thus $Q^\prime \in U$ and hence $Q\in S_r\subseteq U$ which completes the proof.
 \end{proof}
Let $R$ be a $G$-graded ring. As usual, the nilradical of $R$ and the set of all units of $R$ will be denoted by $N(R)$ and $U(R)$, respectively.

\begin{proposition}
Let $M$ be a $G$-graded $R$-module and $r\in h(R)$. Then,
\begin{enumerate}
\item $\rho^{-1}(D_{\overline{r}})=S_r$
\item $\rho(S_r)\subseteq D_{\overline{r}}$. If $\rho$ is surjective, then the equality holds.
\item  $S_r\cap S_t=S_{rt}$, for any $r, t\in h(R)$.
\item If $r\in N(R)$, then $S_r=\emptyset$.
\item  If $r\in U(R)$, then $S_r=\mathcal{PS}_G(M)$.
\end{enumerate}
\end{proposition}
\begin{proof}
(1) $\rho^{-1}(D_{\overline{r}})=\rho^{-1}(Spec_G(\overline{R})-V_G^{\overline{R}}(\overline{r}\overline{R}))=\mathcal{PS}_G(M)-\rho^{-1}(V_G^{\overline{R}}(\overline{r}\overline{R}))=\mathcal{PS}_G(M)-\nu_G(rM)=S_r$ by Proposition 2.10.\\
(2) Trivial.\\
(3) By \cite[Theorem 2.3 (3)]{14}, we have $D_{\overline{r}}\cap D_{\overline{t}}=D_{\overline{rt}}$ and hence $S_r\cap S_t=\rho^{-1}(D_{\overline{r}})\cap \rho^{-1}(D_{\overline{t}})=\rho^{-1}(D_{\overline{rt}})=S_{rt}$.\\
(4)  Assume that $r\in N(R)$. It follows that $D_r=\emptyset$ by \cite[Proposition 3.6 (2)]{18}, and thus $D_{\overline{r}}=\emptyset$. Therefore $S_r=\rho^{-1}(D_{\overline{r}})=\emptyset$ by (1). \\
(5)  Assume that $r\in U(R)$. By\cite[Proposition 3.6 (3)]{18}, we have $D_r=Spec_G(R)$ and hence $D_{\overline{r}}=Spec_G(\overline{R})$. By (1), we obtain $S_r=\rho^{-1}(D_{\overline{r}})=\rho^{-1}(Spec_G(\overline{R}))=\mathcal{PS}_G(M)$.\\
\end{proof}

In part (a) of the next example, we see that if $F$ is a $G$-graded field and $M$ is a $G$-graded $F$-module, then the Zariski topology on $\mathcal{PS}_G(M)$ is the trivial topology. However, if we have a $G$-graded $R$-module $M$ and the Zariski topology on $\mathcal{PS}_G(M)$ is the trivial topology, then $R$ might not be a graded field and this will be discussed in part (b).

\begin{example}(a) Let $F$ be a $G$-graded field and $M$ be a $G$-graded $F$-module. Then any non-zero homogeneous element of $F$ is unit. By Proposition 3.2 (5), we have $S_r=\mathcal{PS}_G(M)$ for any non-zero homogeneous element $r$ of $F$. Also $S_0=\mathcal{PS}_G(M)-\nu_G(0)=\emptyset$ and hence $S=\{S_r\,\mid\,r\in h(F)\}=\{\mathcal{PS}_G(M),\emptyset\}$. Therefore, the Zariski topology on $\mathcal{PS}_G(M)$ is the trivial topology on $\mathcal{PS}_G(M)$. \\
(b) Let $R=\mathbb{Z}_8$ as a $\mathbb{Z}_2$-graded $\mathbb{Z}_8$ module by $R_0=\mathbb{Z}_8$ and $R_1=\{0\}$. Note that $1, 3, 5, 7\in h(\mathbb{Z}_8)\cap U(\mathbb{Z}_8)$. So $S_1=S_3=S_5=S_7=\mathcal{PS}_{\mathbb{Z}_2}(\mathbb{Z}_8)$ by Proposition 3.2 (5). Also $S_0=S_2=S_4=S_6=\emptyset$ by Proposition 3.5 (4), since $0, 2, 4, 6\in N(\mathbb{Z}_8)\cap h(\mathbb{Z}_8)$. Now $S=\{S_r\,\mid\,r\in h(R)\}=\{\emptyset,\, \mathcal{PS}_{\mathbb{Z}_2}(\mathbb{Z}_8)\}$ and hence the Zariski topology on $\mathcal{PS}_{\mathbb{Z}_2}(\mathbb{Z}_8)$ is the trivial topology. But $\mathbb{Z}_8$ is not $\mathbb{Z}_2$-graded field.
\end{example}
\begin{theorem} Let $M$ be a $G$-graded $R$-module. If $\rho$ is surjective, then the open set $S_r$ in $\mathcal{PS}_G(M)$ for each $r\in h(R)$ is quasi compact; in particular, the space $\mathcal{PS}_G(M)$ is quasi compact.
\end{theorem}

\begin{proof} Let $r\in h(R)$ and $\zeta=\{S_t\,\mid\,t\in\Delta\}$ be a basic open cover for $S_r$, where $\Delta$ is an indexing set subset of $h(R)$. Then $S_r\subseteq\underset{t\in\Delta}{\bigcup}S_t$ and thus $D_{\overline{r}}=\rho(S_r)\subseteq\underset{t\in\Delta}{\bigcup}\rho(S_t)=\underset{t\in\Delta}{\bigcup}D_{\overline{t}}$ by Proposition 3.2 (2). Then $\overline{\zeta}=\{D_{\overline{t}}\,\mid\,t\in\Delta\}$ is a basic open cover for the quasi compact set $D_{\overline{r}}$ and hence it has a finite subcover $\overline{\overline{\zeta}}=\{D_{\overline{t_i}}\,\mid\,i=1,...,n\}$, where $t_i\in \Delta$ for any $i=1,...,n$. This means that $D_{\overline{r}}\subseteq\bigcup\limits_{i=1}^n D_{\overline{t_i}}$ and it follows that $S_r=\rho^{-1}(D_{\overline{r}})\subseteq\bigcup\limits_{i=1}^n \rho^{-1}(D_{\overline{t_i}})=\bigcup\limits_{i=1}^n S_{\overline{t_i}}$ by Proposition 3.2 (1). Therefore $\underline{\zeta}=\{S_{\overline{t_i}}\,\mid\,i=1,...,n\}\subseteq\zeta$ is a finite subcover for $S_r$. For the other part of the theorem, since $\mathcal{PS}_G(M)=S_1$, then $\mathcal{PS}_G(M)$ is quasi compact.
\end{proof}
\begin{theorem}
Let $M$ be a $G$-graded $R$-module. If $\rho$ is surjective, then the quasi compact open sets of $\mathcal{PS}_G(M)$ are closed under finite intersection and form an open base.
\end{theorem}
\begin{proof}
Let $C_1, C_2$ be quasi compact open sets of $\mathcal{PS}_G(M)$ and $\zeta=\{S_r\,\mid\,r\in\Delta\}$ be a basic open cover for $C_1\cap C_2$, where $\Delta\subseteq h(R)$ is an indexing set. Since $S$ is a base for the Zariski topology on $\mathcal{PS}_G(M)$, then the quasi compact open sets $C_1, C_2$ can be written as a finite union of elements of $S$. So let $C_1=\bigcup\limits_{i=1}^n S_{t_i}$ and $C_2=\bigcup\limits_{j=1}^m S_{z_j}$. By Proposition 3.2 (3), we have $C_1\cap C_2=\underset{i,j}{\bigcup}(S_{t_i}\cap S_{z_j})=\underset{i,j}{\bigcup}S_{t_i z_j}\subseteq \underset{r\in \Delta}{\bigcup} S_r$. Note that for any $i, j$ we have $t_i z_j\in h(R)$ as $t_i, z_j\in h(R)$. So without loss of generality we can assume that $C_1\cap C_2=\bigcup\limits_{k=1}^L S_{h_k}$ where $h_k\in h(R)$, for $k=1,...,L$. Then $S_{h_k}\subseteq\underset{r\in\Delta}{\bigcup}S_r$  for each $k$. Now each $S_{h_k}$ is quasi compact by Theorem 3.4 and this follows that $S_{h_k}\subseteq \bigcup\limits_{i=1}^{d_k} S_{r_{k,i}}$, where $d_k\geq 1$ depends on $k$ and $r_{k,i}\in \Delta$, for any $k=1,...,L$ and $i=1,...,d_k$. Therefore $C_1 \cap C_2=\bigcup\limits_{k=1}^L S_{h_k}\subseteq\bigcup\limits_{k=1}^L \bigcup\limits_{i=1}^{d_k}S_{r_{k,i}}$ and thus $\overline{\zeta}=\{S_{r_{k,i}}\,\mid\, k=1,...,L, i=1,...,d_k\}$ is a finite subcover for $C_1\cap C_2$. The other part of the theorem is trivially true.
\end{proof}
\section{Irreducibility in $\mathcal{PS}_G(M)$}
Let $M$ be a $G$-graded $R$-module and $Y$ be a subset of $\mathcal{PS}_G(M)$. We will denote the closure of $Y$ in $\mathcal{PS}_G(M)$ by $Cl(Y)$ and the intersection $\underset{Q\in Y}{\bigcap}Gr_M(Q)$ by $\eta(Y)$. If $Z$ is a subset of $Spec_G(R)$ or $Spec_G(M)$, then the intersection of all members of $Z$ will be expressed by $\gamma(Z)$.

\begin{proposition} Let $M$ be a $G$-graded $R$-module and $Y\subseteq \mathcal{PS}_G(M)$. Then $Cl(Y)=\nu_G(\eta(Y))$. Thus, $Y$ is closed in $\mathcal{PS}_G(M)$ if and only if $\nu_G(\eta(Y))=Y$.
\end{proposition}
\begin{proof}
 Let $\nu_G(N)$ be any closed set containing $Y$, where $N\leq_G M$. Note that $Y\subseteq \nu_G(\eta(Y))$, and so it is enough to show that $\nu_G(\eta(Y))\subseteq \nu_G(N)$. So let $Q\in\nu_G(\eta(Y))$. Then $(\eta(Y):_R M)\subseteq (Gr_M(Q):_R M)$. Note that for any $Q^\prime\in Y$, we have $(N:_R M)\subseteq (Gr_M(Q^\prime):_R M)$ and hence $(N:_R M)\subseteq \underset{Q^\prime\in Y}{\bigcap}(Gr_M(Q^\prime):_R M)=(\underset{Q^\prime\in Y}{\bigcap}Gr_M(Q^\prime) :_R M)=(\eta(Y):_R M)\subseteq (Gr_M(Q):_R M)$. Thus $Q\in \nu_G(N)$. Therefore $\nu_G(\eta(Y))$ is the smallest closed set containing $Y$ and hence $Cl(Y)=\nu_G(\eta(Y))$.
 \end{proof}

 Recall that a topological space $X$ is irreducible if any two non-empty open subsets of $X$ intersect. Equivalently, $X$ is irreducible if for any decomposition $X=F_1 \cup F_2$ with closed subsets $F_i$ of $X$ with $i=1, 2$, we have $F_1=X$ or $F_2=X$. A subset $X^\prime$ of $X$ is irreducible if it is an irreducible topological space with the induced topology. Let $X$ be a topological space. Then a subset $Y$ of $X$ is irreducible if and only if its closure is irreducible. Also every singleton subset of $X$ is irreducible, (see \cite{6}).
 \begin{theorem} Let $M$ be a $G$-graded $R$-module. Then for each $Q\in \mathcal{PS}_G(M)$, the closed set $\nu_G(Q)$ is irreducible closed subset of $\mathcal{PS}_G(M)$. In particular, if $\{0\}\in\mathcal{PS}_G(M)$, then $\mathcal{PS}_G(M)$ is irreducible.
 \end{theorem}

\begin{proof}
 For any $Q\in \mathcal{PS}_G(M)$, we have $Cl(\{Q\})=\nu_G(\eta(\{Q\}))=\nu_G(Gr_M(Q))=\nu_G(Q)$ by Proposition 2.5. Now $\{Q\}$ is irreducible in $\mathcal{PS}_G(M)$, then its closure $\nu_G(Q)$ is irreducible. The other part of the theorem follows from the equality $\nu_G(\{0\})=\mathcal{PS}_G(M)$.
 \end{proof}
 In Theorem 4.2, if we drop the condition that $Q\in \mathcal{PS}_G(M)$, then $\nu_G(Q)$ might not be irreducible. Actually, $\mathcal{PS}_G(M)$ itself is not always irreducible. For this, if we take $R=\mathbb{Z}_6$ as $\mathbb{Z}_2$-graded $\mathbb{Z}_6$ module by $R_0=R$ and $R_1=\{0\}$. Then it can easily be checked that $\mathcal{PS}_{\mathbb{Z}_2}(\mathbb{Z}_6)=\{\{0,3\},\{0,2,4\}\}$, $\nu_{\mathbb{Z}_2}(3\mathbb{Z}_6)=\{\{0,3\}\}$ and $\nu_{\mathbb{Z}_2}(2\mathbb{Z}_6)=\{\{0,2,4\}\}$. Now $\mathcal{PS}_{\mathbb{Z}_2}(\mathbb{Z}_6)=\nu_{\mathbb{Z}_2}(3\mathbb{Z}_6)\cup\nu_{\mathbb{Z}_2}(2\mathbb{Z}_6)$. But $\mathcal{PS}_{\mathbb{Z}_2}(\mathbb{Z}_6)\neq\nu_{\mathbb{Z}_2}(3\mathbb{Z}_6)$ and $\mathcal{PS}_{\mathbb{Z}_2}(\mathbb{Z}_6)\neq \nu_{\mathbb{Z}_2}(2\mathbb{Z}_6)$. Therefore $\mathcal{PS}_{\mathbb{Z}_2}(\mathbb{Z}_6)=\nu_{\mathbb{Z}_2}(\{0\})$ is not irreducible.

Now, we need the following lemma to prove the next theorem.

\begin{lemma} \cite[Lemma 4.9]{10} A subset $Y$ of $Spec_G(R)$ for any graded ring $R$ is irreducible if and only if $\gamma(Y)$ is a graded prime ideal of $R$.
\end{lemma}
\begin{proof}
$\Rightarrow$: Let $Y$ be irreducible subset of $Spec_G(R)$ and $r_1, r_2\in h(R)$ with $r_1 r_2\in\gamma(Y)$. Then $r_1 r_2\in p$ for any $p\in Y$. Let $U_1=Y\cap(Spec_G(R)-V_G^R(r_1 R))$ and $U_2=Y\cap(Spec_G(R)-V_G^R(r_2 R))$. If $U_1, U_2$ are non-empty sets, then $U_1\cap U_2\neq\emptyset$ as $Y$ is irreducible and $U_1, U_2$ are open sets in $Y$. So $\exists p\in Y$ such that $r_1 R\nsubseteq p$ and $r_2 R\nsubseteq p$. It follows that $r_1\notin p$ and $r_2\notin p$ and hence $r_1 r_2\notin p$ as $p\in Spec_G(R)$, a contradiction. Therefore $U_1=\emptyset$ or $U_2=\emptyset$. If $U_1=\emptyset$, then $Y\subseteq V_G^R(r_1 R)$. This implies that $r_1 R\subseteq Q$ for any $Q\in Y$ and thus $r_1 R\subseteq\underset{Q\in Y}{\bigcap}Q=\gamma(Y)$. Therefore $r_1\in \gamma(Y)$. Similarly, if $U_2\neq\emptyset$, then $r_2\in\gamma(Y)$. Hence $r_1\in \gamma(Y)$ or $r_2\in \gamma(Y)$.\\
$\Leftarrow$: Assume that $\gamma(Y)$ is a graded prime ideal of $R$, where $Y\subseteq Spec_G(R)$. Let $Y=F_1\cup F_2$, where $F_1, F_2$ are closed sets in $Y$. Now $F_1=V_G^R(I_1)\cap Y$ and $F_2=V_G^R(I_2)\cap Y$ for some $I_1, I_2\lhd_G R$. It follows that $Y=Y\cap V_G^R(I_1\cap I_2)$ and hence $Y\subseteq V_G^R(I_1\cap I_2)$, which implies that $I_1\cap I_2\subseteq p$ for any $p\in Y$. Thus $I_1\cap I_2\subseteq \gamma(Y)$. Since $\gamma(Y)\in Spec_G(R)$, then $I_1\subseteq \gamma(Y)$ or $I_2\subseteq\gamma(Y)$. If $I_1\subseteq\gamma(Y)$, then $V_G^R(\gamma(Y))\subseteq V_G^R(I_1)$ and thus $Y\subseteq V_G^R(I_1)$ as $Y\subseteq V_G^R(\gamma(Y))$. It follows that $F_1=Y$. Similarly, if $I_2\subseteq \gamma(Y)$, we obtain $F_2=Y$. This proves that $F_1=Y$ or $F_2=Y$, hence, $Y$ is irreducible.
\end{proof}

\begin{theorem}
Let $M$ be a $G$-graded $R$-module and $Y\subseteq\mathcal{PS}_G(M)$. Then:
 \begin{enumerate}
\item If $\eta(Y)$ is graded primary submodule of $M$, then $Y$ is irreducible.
\item If $Y$ is irreducible, then $\Upsilon=\{(Gr_M(Q):_R M)\,\mid\,Q\in Y\}$ is an irreducible subset of $Spec_G(R)$, i.e., $\gamma(\Upsilon)=(\eta(Y):_R M)\in Spec_G(R)$.
\end{enumerate}
\end{theorem}
\begin{proof}(1) Assume that $\eta(Y)$ is a graded primary submodule of $M$, then it is easy to see that $\eta(Y)\in \mathcal{PS}_G(M)$. By Theorem 4.2 and Proposition 4.1, we have $\nu_G(\eta(Y))=Cl(Y)$ is irreducible and hence $Y$ is irreducible.\\
(2) Assume that $Y$ is irreducible. Then $\rho(Y)=Y^\prime$ is irreducible subset of $Spec_G(\overline{R})$ as $\rho$ is continuous by Proposition 2.10. Note that $\gamma(Y^\prime)=\gamma(\rho(Y))=\underset{Q\in Y}{\bigcap} \overline{(Gr_M(Q):_R M)}=\overline{(\underset{Q\in Y}{\bigcap}Gr_M(Q):_R M)}=\overline{(\eta(Y):_R M)}$ and hence $\gamma(Y^\prime)=\overline{(\eta(Y):_R M)}\in Spec_G(\overline{R})$ by Lemma 4.3. It follows that $(\eta(Y):_R M)\in Spec_G(R)$. Now $\gamma(\Upsilon)=\underset{Q\in Y}{\bigcap}(Gr_M(Q):_R M)=(\underset{Q\in Y}{\bigcap}Gr_M(Q):_R M)=(\eta(Y):_R M)\in Spec_G(R)$ and thus $\Upsilon$ is irreducible subset of $Spec_G(R)$ by Lemma 4.3 again.
\end{proof}
 Let $X$ be a topological space and $Y$ be a closed subset of $X$. An element $y\in Y$ is called a generic point if $Y=Cl(\{y\})$. An irreducible component of $X$ is a maximal irreducible subset of $X$. The irreducible components of $X$ are closed and they cover $X$, (see \cite{9}).

\begin{theorem}
 Let $M$ be a $G$-graded $R$-module. Let $Y\subseteq \mathcal{PS}_G(M)$ and $\rho$ be surjective. Then $Y$ is an irreducible closed subset of $\mathcal{PS}_G(M)$ if and only if $Y=\nu_G(Q)$ for some $Q\in\mathcal{PS}_G(M)$. Hence every irreducible closed subset of $\mathcal{PS}_G(M)$ has a generic point.
 \end{theorem}
 \begin{proof} Assume that $Y$ is an irreducible closed subset of $\mathcal{PS}_G(M)$. Then $Y=\nu_G(N)$ for some $N\leq_G M$. Also $(\eta(Y):_R M)=(\eta(\nu_G(N)):_R M)\in Spec_G(R)$ by Theorem 4.4. It follows that $\overline{(\eta(Y):_R M)}\in Spec_G(\overline{R})$ and hence $\exists Q\in \mathcal{PS}_G(M)$ such that $(Gr_M(Q):_R M)=(\eta(\nu_G(N)):_R M)$ as $\rho$ is surjective. So $Gr((Q:_R M))=Gr((\eta(\nu_G(N)):_R M))$ and so $\nu_G(Q)=\nu_G(\eta(\nu_G(N)))=Cl(\nu_G(N))=\nu_G(N)=Y$ by Proposition 4.1 and Lemma 2.6 (3). Conversely, if $Y=\nu_G(Q)$ for some $Q\in \mathcal{PS}_G(M)$, then $Y$ is irreducible by Theorem 4.2.
 \end{proof}

\begin{theorem}
Let $M$ be a $G$-graded $R$-module and $Q\in\mathcal{PS}_G(M)$. If $\overline{(Gr_M(Q):_R M)}$ is a minimal graded prime ideal of $\overline{R}$, then $\nu_G(Q)$ is irreducible component of $\mathcal{PS}_G(M)$. The converse is true if $\rho$ is surjective.
\end{theorem}
\begin{proof}
 Note that $\nu_G(Q)$ is irreducible by Theorem 4.2 and it remains to show that it is a maximal irreducible. Let $Y$ be irreducible subset of $\mathcal{PS}_G(M)$ with $\nu_G(Q)\subseteq Y$ and if we show that $Y=\nu_G(Q)$, then we are done. Since $Q\in \nu_G(Q)\subseteq Y$, then $Q\in Y$ and thus $\overline{(\eta(Y):_R M)}\subseteq \overline{(Gr_M(Q):_R M)}$. It follows that $\overline{(\eta(Y):_R M)}=\overline{(Gr_M(Q):_R M)}$ as $\overline{(Gr_M(Q):_R M)}$ is a minimal graded prime ideal of $\overline{R}$ and $\overline{(\eta(Y):_R M)}\in Spec_G(\overline{R})$ by Theorem 4.4 (2). Hence $V_G^{\overline{R}}(\overline{(\eta(Y):_R M)})=V_G^{\overline{R}}(\overline{(Gr_M(Q):_R M)})$, which implies that $\nu_G(\eta(Y))=\rho^{-1}(V_G^{\overline{R}}(\overline{(\eta(Y):_R M)}))=\rho^{-1}(V_G^{\overline{R}}(\overline{(Gr_M(Q):_R M)}))=\nu_G(Gr_M(Q))=\nu_G(Q)$ by Proposition 2.10, Lemma 2.6 (4) and Proposition 2.5. Since $Y\subseteq\nu_G(\eta(Y))$, then $Y\subseteq \nu_G(Q)$ and thus $Y=\nu_G(Q)$. For the converse, we assume that $\rho$ is surjective. Since $Q\in \mathcal{PS}_G(M)$, then $\overline{(Gr_M(Q):_R M)}\in Spec_G(\overline{R})$. Let $\overline{J}\in Spec_G(\overline{R})$ with $\overline{J}\subseteq\overline{(Gr_M(Q):_R M)}$ and it is enough to show that $\overline{J}=\overline{(Gr_M(Q):_R M)}$. Note that $\exists Q^\prime \in \mathcal{PS}_G(M)$ such that $\rho(Q^\prime)=\overline{J}$ as $\rho$ is surjective. So we have $J=(Gr_M(Q^\prime):_R M)$. Now $(Gr_M(Q^\prime):_R M)\subseteq (Gr_M(Q):_R M)$ and thus $\nu_G(Q)\subseteq\nu_G(Q^\prime)$. Since $\nu_G(Q)$ is irreducible component and $\nu_G(Q^\prime)$ is irreducible by Theorem 4.2, then $\nu_G(Q)=\nu_G(Q^\prime)$. By Lemma 2.6 (4), we get $(Gr_M(Q):_R M)=(Gr_M(Q^\prime):_R M)=J$ and thus $\overline{J}=\overline{(Gr_M(Q):_R M)}$.
 \end{proof}
 \begin{corollary} Let $M$ be a $G$-graded $R$-module and $\mathcal{K}=\{Q\in\mathcal{PS}_G(M)\,\mid\,\overline{(Gr_M(Q):_R M)}$ is a minimal graded prime ideal of $\overline{R}\}$. If $\rho$ is surjective, then the following hold:
 \begin{enumerate}
\item $T=\{\nu_G(Q)\,\mid\,Q\in\mathcal{K}\}$ is the set of all irreducible components of $\mathcal{PS}_G(M)$.
\item $\mathcal{PS}_G(M)=\underset{Q\in \mathcal{K}}{\bigcup}{\nu_G(Q)}$.
\item $Spec_G(\overline{R})=\underset{Q\in \mathcal{K}}{\bigcup} V_G^{\overline{R}}(\overline{(Q:_R M)})$.
\item $Spec_G(M)=\underset{Q\in \mathcal{K}}{\bigcup}V_G(Q)$.
\item If $\{0\}\in Spec_G(M)$, then the only irreducible component subset of $\mathcal{PS}_G(M)$ is $\mathcal{PS}_G(M)$ itself.
 \end{enumerate}
\end{corollary}
\begin{proof}
(1) follows from Theorem 4.5 and Theorem 4.6.\\
(2) Since any topological space is the union of its irreducible components, then $\mathcal{PS}_G(M)=\underset{Y\in T}{\bigcup} Y=\underset{Q\in \mathcal{K}}{\bigcup}\nu_G(Q)$. \\
(3) Since $\rho$ is surjective, then $Spec_G(\overline{R})=\rho(\mathcal{PS}_G(M))=\rho(\underset{Q\in\mathcal{K}}{\bigcup}\nu_G(Q))=\underset{Q\in \mathcal{K}}{\bigcup}\rho(\nu_G(Q))=\underset{Q\in \mathcal{K}}{\bigcup}V_G^{\overline{R}}(\overline{(Q:_R M)})$ by Proposition 2.11. \\
(4) $Spec_G(M)=\mathcal{PS}_G(M)\cap Spec_G(M)=(\underset{Q\in \mathcal{K}}{\bigcup}\nu_G(Q))\cap Spec_G(M)=\underset{Q\in \mathcal{K}}{\bigcup}V_G(Q)$ by Lemma 2.6 (1).\\
(5) Assume that $\{0\}\in Spec_G(M)$. Then $\overline{(\{0\}:_R M)}\in Spec_G(\overline{R})$. For any $Q\in \mathcal{K}$, we have $\overline{(\{0\}:_R M)}\subseteq \overline{(Gr_M(Q):_R M)}$ and hence $\overline{(\{0\}:_R M)}=\overline{(Gr_M(Q):_R M)}$ as $\overline{(Gr_M(Q):_R M)}$ is a minimal graded prime ideal of $\overline{R}$. Therefore $Gr((Q:_R M))=Gr((\{0\}:_R M))$ and thus $\nu_G(Q)=\nu_G(\{0\})=\mathcal{PS}_G(M)$ by Lemma 2.6 (3). By (1), the set of all irreducible components of $\mathcal{PS}_G(M)$ is $T=\{\nu_G(Q)\,\mid\,Q\in\mathcal{K}\}=\{\mathcal{PS}_G(M)\}$ which completes the proof.
\end{proof}
\begin{proposition}
 Let $R$ be a $G$-graded principal ideal domain and $M$ be a multiplication graded $R$-module. Let $Y\subseteq \mathcal{PS}_G(M)$. If $\eta(Y)$ is a non-zero graded primary submodule of $M$, then $Y\subseteq \mathcal{PS}_G^p(M)$ for some graded maximal ideal $p$ of $R$.
 \end{proposition}
\begin{proof}
 Clearly, $\eta(Y)=Gr_M(\eta(Y))$. Since $\eta(Y)$ is a graded primary submodule of the graded multiplication module $M$, then $\eta(Y)\in Spec_G(M)$ by \cite[Theorem 13]{15} and hence $(\eta(Y):_R M)\in Spec_G(R)$. If $(\eta(Y):_R M)=\{0\}$, then $\eta(Y)=(\eta(Y):_R M)M=\{0\}$, a contradiction. So $(\eta(Y):_R M)$ is a non-zero graded prime ideal in the graded principle ideal domain $R$ and thus $(\eta(Y):_R M)$ is a graded maximal ideal of $R$. It follows that $\eta(Y)$ is a graded maximal submodule of $M$ as $M$ is a graded multiplication module. Now for any $Q\in Y\subseteq\mathcal{PS}_G(M)$, we have $\eta(Y)\subseteq Gr_M(Q)\neq M$ and thus $\eta(Y)=Gr_M(Q)$. This implies that $(Gr_M(Q):_R M)=(\eta(Y):_R M)$ for any $Q\in Y$. Take $p=(\eta(Y):_R M)$. Therefore $Y\subseteq\mathcal{PS}_G^p(M)$.
 \end{proof}
 A topological space $X$ is called a $T_1$-space if every singleton subset of $X$ is closed. A $G$-graded $R$-module $M$ is called graded finitely generated $R$-module if there are $m_1, m_2,...,m_k\in h(M)$ such that $M=\sum\limits_{i=1}^k Rm_i$.

 \begin{proposition}
Let $M$ be a $G$-graded finitely generated $R$-module. If $\mathcal{PS}_G(M)$ is a $T_1$-space, then $\mathcal{PS}_G(M)=Max_G(M)=Spec_G(M)$, where $Max_G(M)$ is the set of all graded maximal submodule of $M$.
\end{proposition}
\begin{proof}
It is clear that $Max_G(M)\subseteq\mathcal{PS}_G(M)$. Now, let $Q\in \mathcal{PS}_G(M)$. Since $\mathcal{PS}_G(M)$ is a $T_1$-space, then $Cl(\{Q\})=\{Q\}$ and so $\nu_G(Q)=\{Q\}$ by Proposition 4.1 and Proposition 2.5. As $M\neq Q$ is a graded finitely generated module, we obtain $M/Q$ is a non-zero graded finitely generated module and hence $\exists N\leq_G M$ with $Q\subseteq N$ such that $N/Q\in Max_G(M/Q)$ by \cite[Lemma 2.7 (ii)]{4}. Now, it is easy to see that $N\in Max_G(M)\subseteq \mathcal{PS}_G(M)$. Since $(Q:_R M)\subseteq (N:_R M)=(Gr_M(N):_R M)$ and $N\in \mathcal{PS}_G(M)$, then $N\in \nu_G(Q)=\{Q\}$ and thus $N=Q\in Max_G(M)$. Therefore $Max_G(M)=\mathcal{PS}_G(M)$. Now, $Spec_G(M)\subseteq\mathcal{PS}_G(M)=Max_G(M)\subseteq Spec_G(M)$. Hence $Spec_G(M)=\mathcal{PS}_G(M)=Max_G(M)$.
\end{proof}
A topological space is called a $T_0$-space if the closure of any two distinct points are distinct.  A topological space is called spectral space if it is homeomorphic to the prime spectrum of a ring equipped with the Zariski topology. Spectral spaces have been characterized by Hochster\cite[Proposition 4]{9} as the topological spaces $X$ which satisfy the following conditions:\\
(a) $X$ is a $T_0$-space.\\
(b) $X$ is quasi compact.\\
(c) The quasi compact open subsets of $X$ are closed under finite intersection and form an open base.\\
(d) each irreducible closed subset of $X$ has a generic point.

\begin{theorem} Let $M$ be a $G$-graded $R$-module. If $\rho$ is surjective, then $\mathcal{PS}_G(M)$ is spectral space if and only if $\mathcal{PS}_G(M)$ is a $T_0$-space.
\end{theorem}
 \begin{proof}
  By Theorem 3.4, Theorem 3.5 and Theorem 4.5.
  \end{proof}
 \begin{theorem}
 Let $M$ be a $G$-graded $R$-module and $\rho$ be surjective. Then the following statements are equivalent:
  \begin{enumerate}
\item $\mathcal{PS}_G(M)$ is a $T_0$-space.
\item If whenever $\nu_G(Q)=\nu_G(Q^\prime)$ with $Q, Q^\prime\in \mathcal{PS}_G(M)$, then $Q=Q^\prime$.
\item $\rho$ is injective.
\item $|\mathcal{PS}_G^p(M)|\leq 1$ for every $p\in Spec_G(R)$.
\item $\mathcal{PS}_G(M)$ is a spectral space.
 \end{enumerate}
\end{theorem}

\begin{proof} The equivalence of (2), (3) and (4) is proved in Proposition 2.8. Also (1), (5) are equivalent by Theorem 4.10. For (1)$\Rightarrow$(2), assume that $\nu_G(Q)=\nu_G(Q^\prime)$ for $Q, Q^\prime\in \mathcal{PS}_G(M)$, then $Cl(\{Q\})=\nu_G(Q)=\nu_G(Q^\prime)=Cl(\{Q^\prime\})$ and hence $Q=Q^\prime$ as $\mathcal{PS}_G(M)$ is $T_0$ space. For (2)$\Rightarrow$(1), let $Q, Q^\prime\in \mathcal{PS}_G(M)$ with $Q\neq Q^\prime$, then by the assumption (2) we have $\nu_G(Q)\neq \nu_G(Q^\prime)$. Hence $Cl(\{Q\})\neq Cl(\{Q^\prime\})$. Therefore $\mathcal{PS}_G(M)$ is a $T_0$-space.
\end{proof}


\bigskip\bigskip\bigskip\bigskip


\begin{thebibliography} {10}

\bibitem{1} K. Al-Zoubi, The graded primary radical of a graded submodules, An.
Stiint. Univ. Al. I. Cuza Iasi. Mat. (N.S.), 1 (2016), 395-402.

\bibitem{2}K. Al-Zoubi, I. Jaradat, and M. Al-Dolat, On graded P-compactly packed
modules, Open Mathematics, 13 (1) (2015).

\bibitem{3} S. E. Atani, On graded prime submodules, Chiang Mai J.
Sci., 33 (1) (2006), 3-7.

\bibitem{4}S. E. Atani and F. Farzalipour, Notes on the graded prime submodules,
Int. Math. Forum. 1(38) (2006), 1871--1880.

\bibitem{5} S.E. Atani and F. Farzalipour, On graded secondary modules, Turk. J.
Math., 31 (2007), 371-378.

\bibitem{6} N. Bourbaki, Commutative Algebra. Chapter 1--7. Springer-Verlag, Berlin
(1989).

\bibitem{7} A. Y. Darani. Topology on $Spec_{g}(M),$ Buletinul Academiei De Stiinte,
3 (67) (2011), 45-53.

\bibitem{8}J. Escoriza and B. Torrecillas, Multiplication Objects in Commutative
Grothendieck Categories, Comm. in Algebra, 26 (6) (1998), 1867-1883.

\bibitem{9} M. Hochster, Prime ideal structure in commutative rings, Trans. Amer.
Math. Soc., 142 (1969), 43-60.

\bibitem{10} M. Jaradat and K. Al-Zoubi, The Quasi-Zariski topology on the graded
quasi-primary spectrum of a graded module over a graded commutative ring,
(submitted), https://arxiv.org/abs/2106.09519

\bibitem{11} C. Nastasescu and V.F. Oystaeyen, Graded and filtered
rings and modules. Lecture Notes in Mathematics Vol. 758, Springer-Verlag,
Berlin, 2004.

\bibitem{12} C. Nastasescu, F. Van Oystaeyen, Graded ring theory, Mathematical
Library 28, North Holand, Amsterdam, 1982.

\bibitem{13} C. Nastasescu and V. F. Oystaeyen, Methods of Graded Rings, Lecture
Notes in Math., Vol. 1836. Berlin-Heidelberg: Springer-Verlag, 2004.

\bibitem{14} N.A. Ozkirisci , K.H. Oral, U. Tekir, Graded prime spectrum of a graded
module, Iran. J. Sci.Technol., 37A3 (2013), 411-420.

\bibitem{15}K.H. Oral, U. Tekir and A.G. Agargun, On Graded prime and primary
submodules, Turk. J. Math., 35 (2011), 159-167.

\bibitem{16} M. Refai, On properties of $G$-spec($R$), Sci. Math.
Jpn. 53 (2001), no. 3 ,411-415.

\bibitem{17} M. Refai and K. Al-Zoubi, On graded primary ideals, Turk. J. Math., 28
(2004), 217-229.

\bibitem{18}M. Refai, M. Hailat and S. Obiedat, Graded radicals on graded prime
spectra, Far East J. of Math. Sci., part I (2000), 59-73.



 \end{thebibliography}
\end{document}